# On Power Law Scaling Dynamics for Time-fractional Phase Field Models during Coarsening

Lizhen Chen [*], Jia Zhao [†] and Hong Wang [‡]

March 13, 2018

## Abstract

In this paper, we study the phase field models with fractional-order in time. The phase field models have been widely used to study coarsening dynamics of material systems with microstructures. It is known that phase field models are usually derived from energy variation, so that they obey some energy dissipation laws intrinsically. Recently, many works have been published on investigating fractional-order phase field models, but little is known of the corresponding energy dissipation laws. In this paper, we focus on the time-fractional phase field models and report that the effective free energy and roughness obey a universal power law scaling dynamics during coarsening. Mainly, the effective free energy and roughness in the time-fractional phase field models scale by following a similar power law as the integer phase field models, where the power is linearly proportional to the fractional order. This universal scaling law is verified numerically against several phase field models, including the Cahn-Hilliard equations with different variable mobilities and molecular beam epitaxy models. This new finding sheds light on potential applications of time fractional phase field models in studying coarsening dynamics and crystal growths.

## 1 Introduction

Phase field models have been fully exploited to study many material, physical and biology systems [11, 22, 25, 30, 30, 32–34, 36]. One important application of phase field models is to study coarsening dynamics [8, 9, 12, 14, 35], which is a widely observed phenomenon in material systems involving micro-structures. It is characterized by the dissipation of excessive Helmholtz free energy and the growth of a characteristic length scale [8]. Physically speaking, the coarsening dynamics could often be characterized by a power law scaling, mainly if we introduce a characteristic length scale $R(t)$, it grows like a power in time, i.e. $R(t) \approx ct^{1/\alpha}$, where $\alpha$ indicates the order, i.e., $R(t)^\alpha - R(t_0)^\alpha \approx c(t - t_0)$ for some positive constant $c$. The coarsening dynamics of phase field models have already been widely investigated both theoretically and numerically [8]. For instance the Cahn-Hilliard equation with smooth double-well potential and phase-dependent diffusion mobility is studied in [8], where the authors show that the coarsening dynamics of the Cahn-Hilliard equation is relative to the mobility coefficient: for CH equation with constant mobility, the free energy scales as $O(t^{-\frac{1}{3}})$; for the two-sided degenerate mobility, the free energy scales $O(t^{-\frac{1}{4}})$, which is independent of the volume fraction of each phase [8]; for the one-sided degenerate mobility, the free energy scales depending on the volume fraction of the phases. The coarsening of the molecular beam epitaxy (MBE) model has also been investigated [14, 15, 29]: for the MBE model with slope selection, the free energy scales as $O(t^{-\frac{1}{3}})$ and the roughness scales as $O(t^{\frac{1}{3}})$; for the MBE model without slope selection, the roughness growth scales as $O(t^{\frac{1}{2}})$.

[*]Beijing Computational Science Research Center, Beijing, China; Email: lzchen@csrc.ac.cn.
[†]Utah State University, Logan, UT, USA; email: jia.zhao@usu.edu.
[‡]University of South Carolina, Columbia, SC, USA; email: hwang@math.sc.edu.

Fractional phase field models have been studied in the last few years and have quickly attracted significantly increased attentions ever since. The study of fractional phase-field models is motivated by the significantly increased applications of fractional partial differential equations (FPDEs) to a wide variety of areas [19, 20, 23] and the rapid developed mathematical, numerical and computational analysis of FPDEs over the past few decades. The concept of FPDEs can be explained in the context of anomalously diffusive transport. The classical second-order diffusion PDE was first presented by Fick in his study on how water and nutrients travel through cell membranes. Einstein and Pearson derived the diffusion PDE independently from first principles and random walk under the common assumptions of (i) the existence of a mean free path and (ii) the existence of a mean time taken to perform a jump in the particle movements in the underlying processes. Under these assumptions, the variance of a particle excursion distance is finite. The central limit theorem concludes that the probability of finding a particle somewhere in space satisfies a Gaussian distribution, which gives a probabilistic description of the Fickian diffusion. In the last few decades more and more diffusion processes were found to be non-Fickian, ranging from the signaling of biological cells [21], anomalous electrodiffusion in nerve cells [13], foraging behavior of animals [24] to viscoelastic and viscoplastic flow [23] and solute transport in groundwater [4]. In these cases either the particle movements have many long jumps or have experienced longer waiting times, namely, the assumptions (i) or (ii) is violated. Consequently, the processes can have large deviations from the stochastic process of Brownian motion, leading to anomalously diffusive transport that exhibits heavy tails either in space or in time. Consequently, classical integer-order PDE models fail to provide an accurate description of these problems. In contrast, it has been shown that these anomalously diffusive transport can be better modeled by space-fractional PDE (for superdiffusive transport) or time-fractional PDE (for subdiffusive transport) [19, 20].

The same principle applies to phase-field models. If the coarsening dynamic processes experience some long-range spatial interaction and so the assumption (i) is violated, the corresponding phase-field models may be better described by space-fractional phase-field models [28]. As a matter of fact, the original phase-field model was expressed in terms of an integration of Helmholtz free energy that has a nonlocal spatial interaction [5]. It was approximated by a Laplacian operator just for the modeling, mathematical and numerical simplicity. On the other hand, when the material during the coarsening process exhibits memory effect, the assumption (ii) is violated. The corresponding phase-field model may be described by a time-fractional phase-field model [17]. In this paper we study time-fractional phase-field models.

In this paper, we lay out a foundation for time fractional phase field models by discovering a universal power laws scaling property. *The scaling of effective free energy/roughness in the time-fractional phase field models during coarsening follow a similar power law as the integer phase field models, where the power is linearly proportional to the fractional order $\alpha$.* By far, there is the first report in literature. Several examples of fractional phase field models with specific choices of free energies are presented to verify this new finding. A rigorous asymptotic analysis to this interesting correlation will be pursed in our future work.

In the rest of the paper, we will introduce the classical and fractional phase field models in Section 2. And in Section 3, we will provide the numerical approximations for the fractional phase field models. Afterward, we will present the power scaling laws for the effective free energy and roughness of the time-fractional phase field models. In the end, we draw a brief conclusion and point out possible future research directions.

## 2 Phase Field Models

In this section, we first briefly introduce the classical (integer) phase field models and the analogical time-fractional phase field models, subsequently. Interested readers could recommended to read [1, 2, 6, 23, 28] and the references therein for more details.

### 2.1 Classical Phase Field Models

Given a domain $\Omega$ with smooth boundary $\partial\Omega$, we consider a material system with two components. We introduce a phase variable $\phi$ to represent the volume fraction of one component, such that $1 - \phi$ represents the volume fraction of the other component. For the situations of material system with multiple components, more phase variables should be introduced. We also introduce the notation $F$ to represent the total free energy. Given the explicit expression of $F$, there are usually two types of phase field equations, namely the Allen-Cahn type (AC) equation and the Cahn-Hilliard type (CH) equation.

The Allen-Cahn equation could be considered as a gradient flow equation, i.e. the rate of change $\partial_t \phi$ is in the direction of decreasing gradient of the free energy functional $F(\phi)$, which reads as

$$\begin{cases} \partial_t \phi = -\lambda \frac{\delta F}{\delta \phi}, \text{ in } \Omega_T = \Omega \times (0, T], \\ \nabla \phi \cdot \mathbf{n} = 0, \text{on } \partial\Omega \times (0, T], \\ \phi = \phi_0, \text{ in } \Omega_0 = \Omega \times \{0\}, \end{cases} \quad (2.1)$$

where $\lambda$ is the motility parameter, and $\frac{\delta F}{\delta \phi}$ is the chemical potential ( the functional deriva-tive of $F$ with respect to $\phi$). Its energy dissipation law is

$$\frac{dF}{dt} = \int_\Omega \frac{\delta F}{\delta \phi}\frac{\delta \phi}{\delta t}d\mathbf{x} = \int_\Omega -\lambda\Big(\frac{\delta F}{\delta \phi}\Big)^2 d\mathbf{x} \leq 0. \tag{2.2}$$

One other property of Allen-Cahn equation is its lacking of mass conservation, i.e.

$$\frac{d}{dt}\int_\Omega \phi d\mathbf{x} = \int_\Omega -\lambda\frac{\delta F}{\delta \phi}d\mathbf{x} \neq 0. \tag{2.3}$$

The Cahn-Hilliard equation could be interpreted as derived by following the Fick's first law, where the flux goes from regions of high chemical potential to low chemical potential, which reads as

$$\begin{cases} \partial_t \phi = \nabla \cdot (\lambda \nabla \frac{\delta F}{\delta \phi}), \text{in } \Omega_T = \Omega \times (0,T], \\ \nabla \phi \cdot \mathbf{n} = 0, \nabla \frac{\delta F}{\delta \phi} \cdot \mathbf{n} = 0, \text{on } \partial\Omega \times (0,T], \\ \phi = \phi_0, \text{in } \Omega_0 = \Omega \times \{0\}, \end{cases} \tag{2.4}$$

where $\lambda$ is the motility parameter. For the Cahn-Hilliard equation, it possesses energy dissipation property as follows

$$\frac{dF}{dt} = \int_\Omega \frac{\delta F}{\delta \phi}\frac{\delta \phi}{\delta t}d\mathbf{x} = -\int_\Omega \lambda|\nabla\frac{\delta F}{\delta \phi}|^2 d\mathbf{x} + \int_{\partial\Omega} \lambda\frac{\delta F}{\delta \phi}\nabla\frac{\delta F}{\delta \phi}\cdot \mathbf{n}dS \leq 0, \tag{2.5}$$

and mass conservation property as

$$\frac{d}{dt}\int_\Omega \phi d\mathbf{x} = \int_\Omega \nabla\cdot(\lambda\nabla\frac{\delta F}{\delta \phi})d\mathbf{x} = \int_{\partial\Omega}\lambda\nabla\frac{\delta F}{\delta \phi}\cdot \mathbf{n}ds = 0. \tag{2.6}$$

### 2.2 Time Fractional Phase Field Models

Next, we introduce the corresponding time-fractional phase field models, i.e., the time-fractional Allen-Cahn (FAC) type equations and the time-fractional Cahn-Hilliard (FCH) type equations. Mainly, we simply replace the integer time derivatives by the fractional time derivatives.

For simplicity, we consider periodic boundary conditions for the FPDEs in this paper. Thus, the fractional Allen-Cahn equation reads

$$\begin{cases} {}^C_0D_t^\alpha \phi = -\lambda\frac{\delta F}{\delta \phi}, \text{in } \Omega_T = \Omega \times (0,T], \\ \phi = \phi_0, \text{in } \Omega_0 = \Omega \times \{0\}. \end{cases} \tag{2.7}$$

and the fractional Cahn-Hilliard equation reads

$$\begin{cases} {}_0^C D_t^\alpha \phi = \nabla \cdot (\lambda \nabla \frac{\delta F}{\delta \phi}), \text{in } \Omega_T = \Omega \times (0, T], \\ \phi = \phi_0, \text{in } \Omega_0 = \Omega \times \{0\}. \end{cases} \quad (2.8)$$

Here ${}_0^C D_t^\alpha$ (with $\alpha \in (0,1]$) is the Caputo fractional derivative of order $\alpha$ defined by [16, 23]

$$ {}_0^C D_t^\alpha \phi(t) := \begin{cases} \dfrac{1}{\Gamma(1-\alpha)} \displaystyle\int_0^t \dfrac{\phi_\eta(\eta)}{(t-\eta)^\alpha} d\eta, & 0 \leq \alpha < 1, \\ \dfrac{\partial \phi}{\partial t}, & \alpha = 1, \end{cases} \quad (2.9)$$

with $\Gamma(\cdot)$ the Gamma function. The Caputo type time-fractional differential operator provides a means to model the sub-diffusive or long time memory behavior.

*Remark* 2.1. There usually exist two ways to express the operator ${}_0 D_t^\alpha$ : Caputo and Riemann-Liouville definitions. In the Caputo definition, the fractional derivative of order $\alpha$, denoted by ${}_0^C D_t^\alpha u(x,t)$ is defined as

$$ {}_0^C D_t^\alpha u(x,t) = \frac{1}{\Gamma(1-\alpha)} \int_0^t \frac{\partial u(x,\tau)}{\partial \tau} \frac{d\tau}{(t-\tau)^\alpha}, \quad 0 < \alpha < 1, \quad (2.10)$$

where $\Gamma(\cdot)$ denotes Gamma function. By using integration by parts, we have the alternative formula

$$ {}_0^C D_t^\alpha u(x,t) = \frac{1}{\Gamma(1-\alpha)} \left( \frac{u(x,t) - u(x,0)}{t^\alpha} + \alpha \int_0^t \frac{u(x,t) - u(x,s)}{(t-s)^{\alpha-1}} ds \right). \quad (2.11)$$

The Riemann-Liouville derivative ${}_0^R D_t^\alpha$ has the form

$$ {}_0^R D_t^\alpha u(x,t) = \frac{1}{\Gamma(1-\alpha)} \frac{\partial}{\partial t} \int_0^t \frac{u(x,\tau)}{(t-\tau)^\alpha} d\tau, \quad 0 < \alpha \leq 1. \quad (2.12)$$

The two definitions are linked by

$$ {}_0^R D_t^\alpha u(x,t) = \frac{u(x,0)}{\Gamma(1-\alpha)t^\alpha} + {}_0^C D_t^\alpha u(x,t). \quad (2.13)$$

In this paper, we stick to the Caputo's definition for simplicity.

## 3 Accurate and Efficient Numerical Approximation

In order to solve the time-fractional phase field models, accurate, efficient and stable numerical approximations are necessary. There are mainly two major difficulties in dis-cretizing the models in time, mainly dealing with the time-fractional operator, and dealing with the nonlinear terms on the right-hand-side of the models. Here we provide general numerical techniques to deal with these issues as follows.

Recall the Caputo fractional derivative of order $\alpha$ $(0 < \alpha < 1)$ [3, 18]

$$_0^C D_t^\alpha u(x,t) = \frac{1}{\Gamma(1-\alpha)} \int_0^t \frac{\partial u(x,\eta)}{\partial \eta}(t-\eta)^{-\alpha} d\eta. \tag{3.1}$$

We consider the fast evaluation of the Caputo fractional derivative proposed by [10, 31]. Before introduce the scheme, we approximated the kernel $t^{-\beta}(0 < \beta < 2)$ via a sum-of-exponentials approximation efficiently on the interval $[\delta, T]$ with $\delta = \min_{1 \leq K \leq N} \Delta t_n$ and the absolute error $\epsilon$. That is to say, there exits positive real numbers $s_i$ and $\omega_i (i = 1, \ldots, K_{exp})$ such that

$$\Big|\frac{1}{t^\beta} - \sum_{i=1}^{K_{exp}} \omega_i e^{-s_i t}\Big| \leq \epsilon, t \in [\delta, T].$$

Then we split the fractional-in-time derivative term into a sum of local part and history part as following:

$$\begin{aligned}
_0^C D_t^\alpha \phi(\mathbf{x}, t_{k+1}) &= \frac{1}{\Gamma(1-\alpha)} \int_0^t \frac{\phi_s(\mathbf{x},s)}{(t-s)^\alpha} ds, \\
&= \frac{1}{\Gamma(1-\alpha)} \int_{t_k}^{t_{k+1}} \frac{\phi_s(\mathbf{x},s)}{(t-s)^\alpha} ds + \frac{1}{\Gamma(1-\alpha)} \int_0^{t_k} \frac{\phi_s(\mathbf{x},s)}{(t-s)^\alpha} ds \\
&:= C_l(t_{k+1}) + C_h(t_{k+1}),
\end{aligned} \tag{3.2}$$

The local part $C_l(t_{k+1})$ and the history part $C_h(t_{k+1})$ are approximated respectively as follows:

$$C_l(t_{k+1}) \approx \frac{\phi(\mathbf{x}, t_{k+1}) - \phi(\mathbf{x}, t_k)}{\Delta t_{k+1}^\alpha \Gamma(2-\alpha)},$$

and

$$\begin{aligned}
C_h(t_{k+1}) &= \frac{1}{\Gamma(1-\alpha)} \Big[\frac{\phi(\mathbf{x}, t_k)}{\Delta t_{k+1}^\alpha} - \frac{\phi(\mathbf{x}, t_0)}{t_{k+1}^\alpha} - \alpha \int_0^{t_k} \frac{\phi(\mathbf{x},s)}{(t-s)^{1+\alpha}} ds\Big] \\
&\approx \frac{1}{\Gamma(1-\alpha)} \Big[\frac{\phi(\mathbf{x}, t_k)}{\Delta t_{k+1}^\alpha} - \frac{\phi(\mathbf{x}, t_0)}{t_{k+1}^\alpha} - \alpha \sum_{i=1}^{K_{exp}} \omega_i \int_0^{t_k} e^{-(t_{k+1}-\tau)s_i} \phi(\mathbf{x}, \tau) d\tau\Big] \\
&= \frac{1}{\Gamma(1-\alpha)} \Big[\frac{\phi(\mathbf{x}, t_k)}{\Delta t_{k+1}^\alpha} - \frac{\phi(\mathbf{x}, t_0)}{t_{k+1}^\alpha} - \alpha \sum_{i=1}^{K_{exp}} \omega_i U_{hist,i}(t_{k+1})\Big],
\end{aligned}$$

with

$$U_{hist,i}(t_{k+1}) = \int_0^{t_k} e^{-(t_{k+1}-\tau)s_i}\phi(\mathbf{x},\tau)d\tau$$

$$= e^{-s_i\Delta t_{k+1}}U_{hist,i}(t_k) + \int_{t_{k-1}}^{t_k} e^{-(t_{k+1}-\tau)s_i}\phi(\mathbf{x},\tau)d\tau$$

$$\approx e^{-s_i\Delta t_{k+1}}U_{hist,i}(t_k) + \frac{e^{-s_i\Delta t_{k+1}}}{s_i^2 \Delta t_k}\Big[(e^{-s_i\Delta t_k} - 1 + s_i\Delta t_k)\phi(\mathbf{x},t_k)$$

$$+ (1 - e^{-s_i\Delta t_k} - e^{-s_i\Delta t_k}s_i\Delta t_k)\phi(\mathbf{x},t_{k-1})\Big].$$

Denote

$$_0^C D_t^\alpha \phi^{n+1} = \frac{\phi^{n+1} - \phi^n}{\Delta t_{n+1}^\alpha \Gamma(2-\alpha)} + \frac{1}{\Gamma(1-\alpha)}\Big[\frac{\phi^n}{\Delta t_{n+1}^\alpha} - \frac{\phi^0}{t_{n+1}^\alpha} - \alpha \sum_{i=1}^{K_{exp}} \omega_i U_{hist,i}^{n+1}\Big],$$

where

$$U_{hist,i}^{n+1} = e^{-s_i\Delta t_{n+1}}U_{hist,i}^n + \frac{e^{-s_i\Delta t_{n+1}}}{s_i^2 \Delta t_n}\Big[(e^{-s_i\Delta t_n} - 1 + s_i\Delta t_n)\phi^n$$

$$+ (1 - e^{-s_i\Delta t_n} - e^{-s_i\Delta t_n}s_i\Delta t_n)\phi^{n-1}\Big].$$

We further utilize the stabilized semi-implicit strategy [26, 27] to discretetize the nonlinear terms on the right hand side to end up with linear schemes. Given the free energy $F$, if we rewrite it as

$$F = (\frac{1}{2}L\phi, \phi) + E, \tag{3.3}$$

where $L$ is the linear operators that are separated from $F$, and $E$ is the rest term.

The scheme for the fractional Allen-Cahn equation is given as follows.

*Scheme* 3.1. Give the initial condition $\phi^0 = \phi_0$. After we obtain $\phi^n$, $n \geq 1$, we can get $\phi^{n+1}$ via

$$C_l(t_{n+1}) + C_h(t_{n+1}) = \lambda\Big[L\phi^{n+1} + S\Delta(\phi^{n+1} - \overline{\phi}^{n+1}) - \overline{\Big(\frac{\delta E}{\delta \phi}\Big)}^{n+1}\Big], \tag{3.4}$$

where $\overline{(\bullet)}^{n+1} = 2(\bullet)^n - (\bullet)^{n-1}$ is a second-order extrapolation, and $S$ is stabilizing constant. The operators $C_l$, $C_h$ are defined in (3.2).

The scheme for the fractional Cahn-Hilliard equation is given as follows.

*Scheme* 3.2. Give the initial condition $\phi^0 = \phi_0$. After we obtain $\phi^n$, $n \geq 1$, we can get $\phi^{n+1}$ via

$$C_l(t_{n+1}) + C_h(t_{n+1}) = \nabla \cdot \Big[\lambda(\overline{\phi}^{n+1})\nabla\Big(L\phi^{n+1} + \overline{\frac{\delta E}{\delta \phi}}^{n+1}\Big)\Big]$$
$$+ \lambda\Big[-S_0\Delta^2(\phi^{n+1} - \overline{\phi}^{n+1}) + S_1\Delta(\phi^{n+1} - \overline{\phi}^{n+1})\Big], \tag{3.5}$$

where $\overline{(\bullet)}^{n+1} = 2(\bullet)^n - (\bullet)^{n-1}$ is a second-order extrapolation, and $S_0$, $S_1$ are stabilizing constant. The operators $C_l$, $C_h$ are defined in (3.2).

## 4 Power Law Scaling Dynamics

In this section, we present our numerical findings on the fractional-order phase field models, that is *the scaling of effective free energy/roughness in the time-fractional phase field models during coarsening follow a similar power law as the integer phase field models, where the power is linearly proportional to the fractional order $\alpha$.*

Here we chose periodic boundary conditions in the square domain. Define the roughness measure function $W(t)$ as follows:

$$W(t) = \sqrt{\frac{1}{|\Omega|} \int_\Omega \left(\phi(x,y,t) - \overline{\phi}(x,y,t)\right)^2 d\Omega}, \tag{4.1}$$

where $\overline{\phi}(x,y,t) = \int_\Omega \phi(x,y,t) d\Omega$. By the linear least square fitting, we get the absolute value of slope for each linear line of energy and roughness, $\beta(\alpha), R(\alpha)$, defined as

$$\log_{10} E(\alpha, t) = \beta^0(\alpha) - \beta(\alpha) \log_{10} t, \quad \log_{10} W(\alpha, t) = R^0(\alpha) + R(\alpha) \log_{10} t. \tag{4.2}$$

**Fractional Cahn-Hilliard equations with the double-well potential.** In the first example, we study the fractional Cahn-Hilliard (FCH) equation, equipped with the free energy

$$F = \int_\Omega \left[\frac{\varepsilon^2}{2} |\nabla \phi|^2 + \frac{1}{4}\phi^2(1-\phi)^2\right] d\mathbf{x}, \tag{4.3}$$

where the first term represents the enthalpy, and the second term is the bulk energy functional. The fractional Cahn-Hilliard equation with a variable mobility reads

$$\begin{cases} {}_0^C D_t^\alpha \phi = \nabla \cdot (\lambda(\phi) \nabla \mu), & (\mathbf{x}, t) \in \Omega \times (0, T], \\ \mu = -\varepsilon^2 \Delta \phi + \phi(1-\phi)(\frac{1}{2} - \phi), & (\mathbf{x}, t) \in \Omega \times (0, T], \\ \phi(\mathbf{x}, 0) = \phi_0(\mathbf{x}), & \mathbf{x} \in \Omega, \end{cases} \tag{4.4}$$

where $\lambda(\phi)$ is the motility parameter, which could be chosen as [8]

$$(i) \lambda(\phi) = \lambda_0, \quad (ii) \lambda(\phi) = \lambda_0 |1 - \phi^2|; \quad (iii) \lambda(\phi) = \frac{\lambda_0}{2}|1 + \phi|. \tag{4.5}$$

First of all, we focus on the FCH equation (4.4) with a constant mobility $\lambda(\phi) = \lambda_0$. We follow the parameters used in [8], i.e. consider domain $[0, 4\pi]^2$, and $\varepsilon = 0.05$, $\lambda_0 = 0.02$, and we choose random initial condition $\phi_{t=0} = 0.001 \times rand(-1, 1)$. The results are summarized in Figure 4.1, where we observe the energy scales approximately like $O(t^{-\frac{\alpha}{3}})$, which is consistent with $O(t^{-\frac{1}{3}})$ as indicated in [8] for the integer Cahn-Hilliard equation with constant mobility. In addition, the coarsening dynamics for the FCH model (2.8) of constant mobility with fractional-order $\alpha = 0.9$ is shown in Figure 4.2.

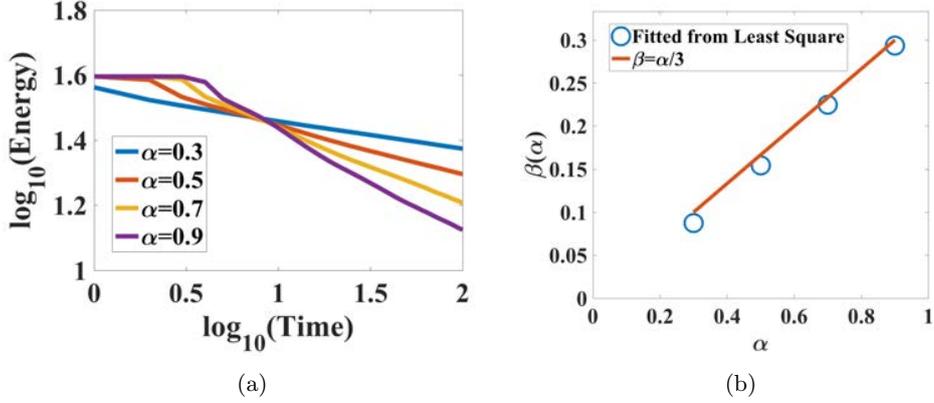

Figure 4.1: Energy scaling laws for the FCH model with constant mobility $\lambda = \lambda_0$. (A) The log-log plot of the effective energy and time with different fractional-order $\alpha$. (b) The least-square fitted energy decay power $\beta(\alpha)$ for different fractional order $\alpha$.

Next, we study the FCH equation (4.4) with variable mobility $\lambda(\phi) = \lambda_0|1-\phi^2|$. We use the same initial values and parameters as the example above. The energy plot at different time slots with various fractional orders are summarized in Figure 4.3. We observe that the energy scales approximately like $O(t^{-\frac{\alpha}{4}})$, which is consistent with $O(t^{-\frac{1}{4}})$ as indicated in [8] for the integer Cahn-Hilliard equation with variable mobility $\lambda(\phi) = \lambda_0|1-\phi^2|$. In addition, the corresponding coarsening dynamics with $\alpha = 0.9$ at different time steps are summarized in Figure 4.4.

**Time-fractional molecular beam epitaxy models.** In the second example, we consider the fractional molecular beam epitaxy (FMBE) model. Given a smooth domain $\Omega$, and use $\phi(\mathbf{x},t) : \Omega \to \mathbb{R}$ to denote the height function of MBE, and the effective free energy is given as

$$E(\phi) = \int_\Omega \left[\frac{\varepsilon^2}{2}|\Delta\phi|^2 + f(\nabla\phi)\right]d\Omega. \qquad (4.6)$$

Here the first term represents the isotropic surface diffusion effect with $\varepsilon$ a constant controlling the surface diffusion strength, and the second term approximates the Enrlich-Schwoebel effect that the adatoms stick to the boundary from an upper terrace, contributing to the steepening of mounds in the film [7].

If we choose $f(\nabla\phi) = -\frac{1}{2}\ln(1 + |\nabla\phi|^2)$, the corresponding FMBE model with slope

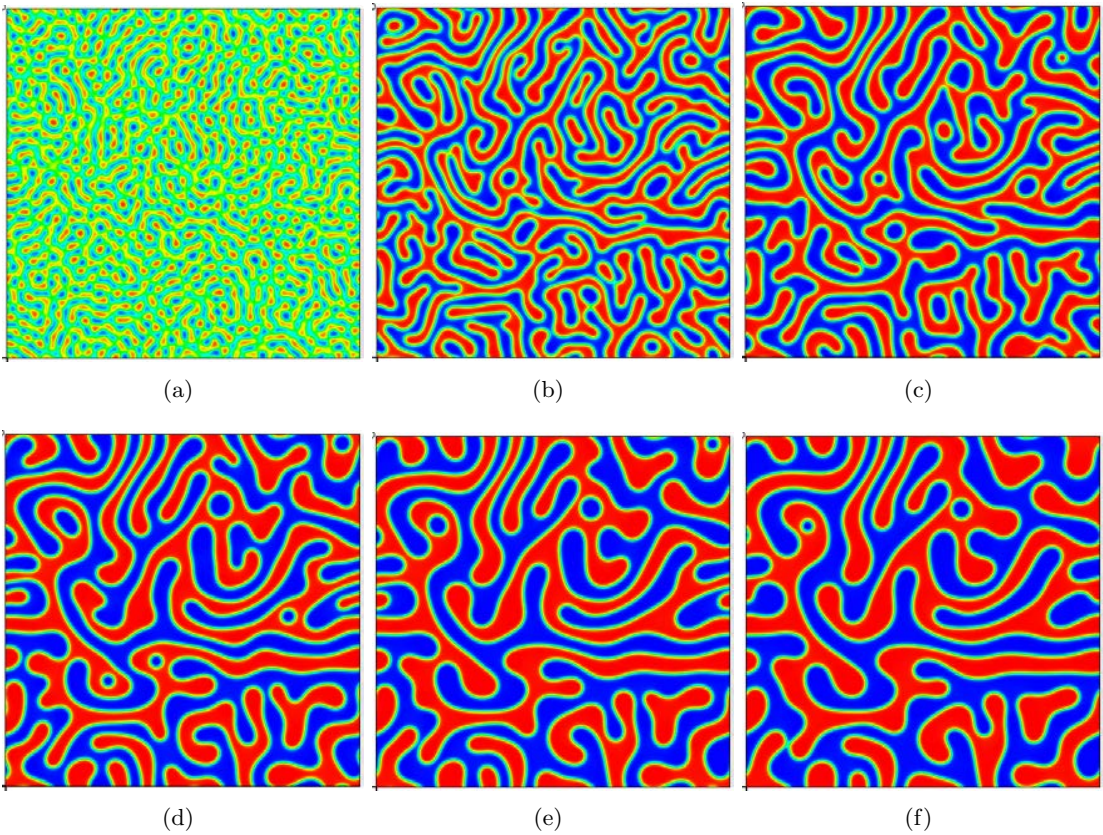

Figure 4.2: Time snapshots of coarsening dynamics driven by the FCH model with constant mobility $\lambda = \lambda_0$. The profiles of $\phi$ at different time slots $t = 5, 25, 50, 100, 125, 145$ are presented.

selection reads

$$\begin{cases} {}^{C}_{0}D^{\alpha}_{t}\phi(\mathbf{x},t) = -M\Big(\varepsilon^2\Delta^2\phi - \nabla\cdot\big((|\nabla\phi|^2 - 1)\nabla\phi\big)\Big), & (\mathbf{x},t) \in \Omega \times (0,T], \\ \phi(\mathbf{x},0) = \phi_0(\mathbf{x}), & \mathbf{x} \in \Omega. \end{cases} \quad (4.7)$$

For the FMBE model with slope selection, the energy scaling laws are summarized in 4.5. It indicates that the energy scales approximately like $O(t^{-\frac{\alpha}{3}})$, which is consistent with $O(t^{-\frac{1}{3}})$ as indicated in [29] for the integer MBE equation with slope selection. Similarly, as shown in Figure 4.6, the roughness scales approximately like $O(t^{\frac{\alpha}{3}})$, which is consistent with $O(t^{\frac{1}{3}})$ as indicated in [29] for the integer MBE equation with slope selection. One time sequence for the time-fractional MBE model with $\alpha = 0.9$ is summarized in Figure 4.7 and 4.8.

On the other hand, if we choose choose $f(\nabla\phi) = \frac{1}{4}(|\nabla\phi|^2 - 1)^2$, the corresponding

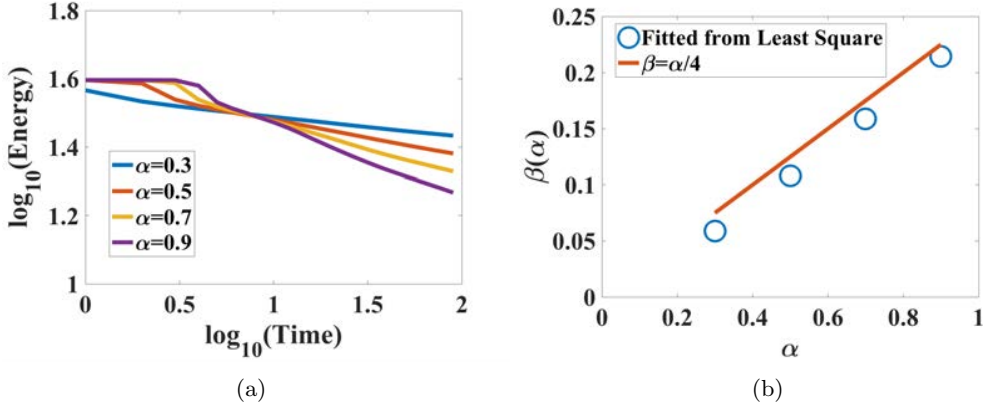

Figure 4.3: Energy Scaling Laws for the FCH mdoel with variable mobility $\lambda = \lambda_0|1 - \phi^2|$. (A) The log-log plot of the effective energy and time with different fractional-order $\alpha$. (b) The least-square fitted energy decay power $\beta$ for different fractional order $\alpha$. This figure illustrates that the energy scales approximately like $O(t^{-\frac{\alpha}{4}})$, which is consistent with $O(t^{-\frac{1}{4}})$ as indicated in [8] for the integer Cahn-Hilliard equation.

time-fractional MBE model without slope selection reads

$$\begin{cases} {}^C_0 D^\alpha_t \phi(\mathbf{x},t) = -M\Big(\varepsilon^2 \Delta^2 \phi + \nabla \cdot \big(\frac{\nabla \phi}{1 + |\nabla \phi|^2}\big)\Big), \ (\mathbf{x},t) \in \Omega \times (0,T], \\ \phi(\mathbf{x},0) = \phi_0(\mathbf{x}), \quad \mathbf{x} \in \Omega. \end{cases} \qquad (4.8)$$

For the time-fractional MBE model without slope selection, the roughness scales approximately like $O(t^{\frac{\alpha}{2}})$, which is consistent with $O(t^{\frac{1}{2}})$ as indicated in [29] for the integer MBE equation without slope selection. The results are summarized in 4.9. One time sequence for the time-fractional MBE model without slope selection of time-fractional order $\alpha = 0.9$ is summarized in Figure 4.10 and 4.11.

## 5 Conclusion

In this paper, we propose time-fractional phase field models, develop their efficient, accurate, full discrete, linear numerical approximation. Correlations of scaling laws with fractional order is studied. The numerical scheme utilizes the fast algorithm for the Caputo fractional derivative operator in time discretization and Fourier spectral method in spatial discretization. Our proposed scheme can handle long time simulation in high dimensions. The numerical approximation strategy proposed in this paper can be readily applied to study many classes of time-fractional and high dimensional phase field models.

With the proposed time fractional Cahn-Hilliard equation, we study the power law scaling behavior of effective free energies with both constant mobility and variable mobility in coarsening dynamics. For time-fractional Cahn-Hilliard equation with constant mobility, it scales as $O(t^{-\frac{\alpha}{3}})$, which is consistent with the well known result $(t^{-\frac{1}{3}})$ for integer order Cahn-Hilliard model [8]. For time fractional Cahn-Hilliard equation with phase-variable dependent mobility $\lambda = \lambda_0|1-\phi^2|$, the effective energy scales as $O(t^{-\frac{\alpha}{4}})$, which is consistent with the result of $O(t^{-\frac{1}{4}})$ for the integer Cahn-Hilliard equation with such variable mobility.

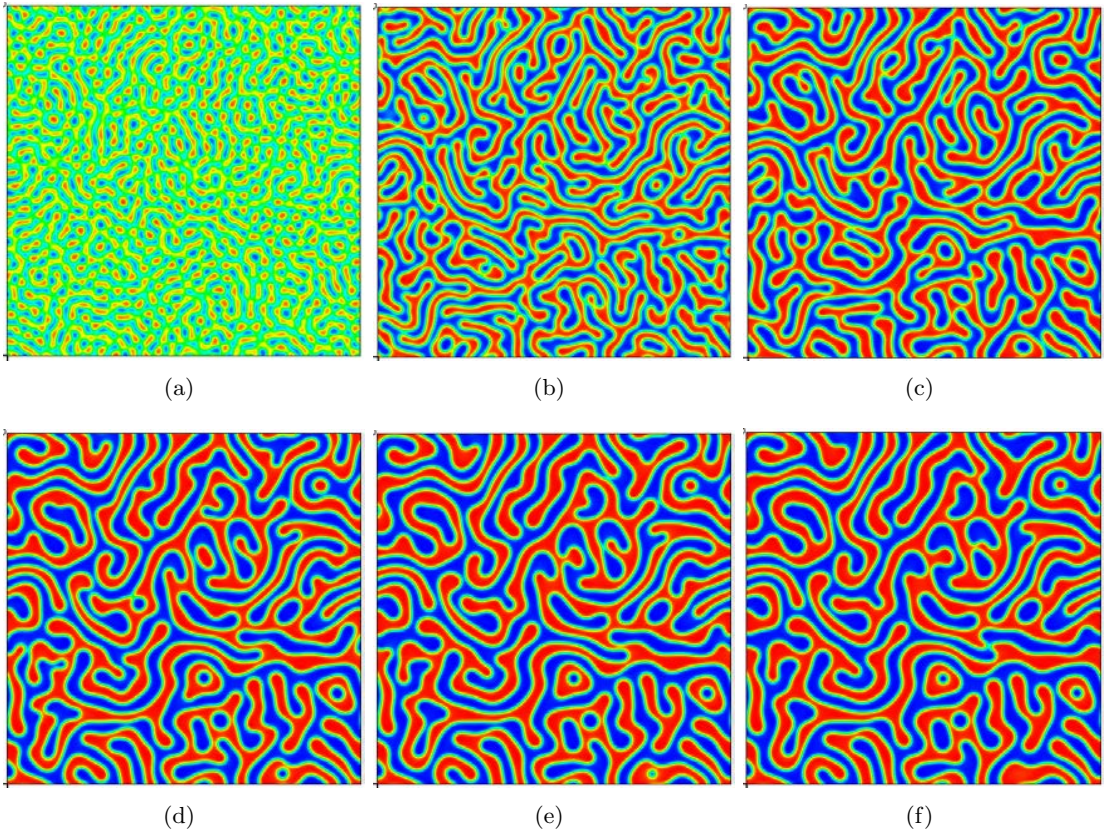

Figure 4.4: Time snapshots of coarsening dynamics driven by the FCH model with variable mobility $\lambda = \lambda_0|1-\phi^2|$. The profiles of $\phi$ at different time slots $t = 5, 25, 50, 100, 125, 145$ are presented.

With the proposed fractional MBE model, we observe that the scaling law for the energy decays as $O(t^{-\frac{\alpha}{3}})$ and the roughness increases as $O(t^{\frac{\alpha}{3}})$ for the MBE model with slope selection, and the roughness increases as $O(t^{\frac{\alpha}{2}})$ for the MBE model without slope selection. That is to say, the coarsening rate of MBE model could be manipulated by the fractional order $\alpha$, and it is linearly proportional to $\alpha$. This is the first time in literature to report/discover such scaling correlation. It provides a potential application field for fractional differential equations.

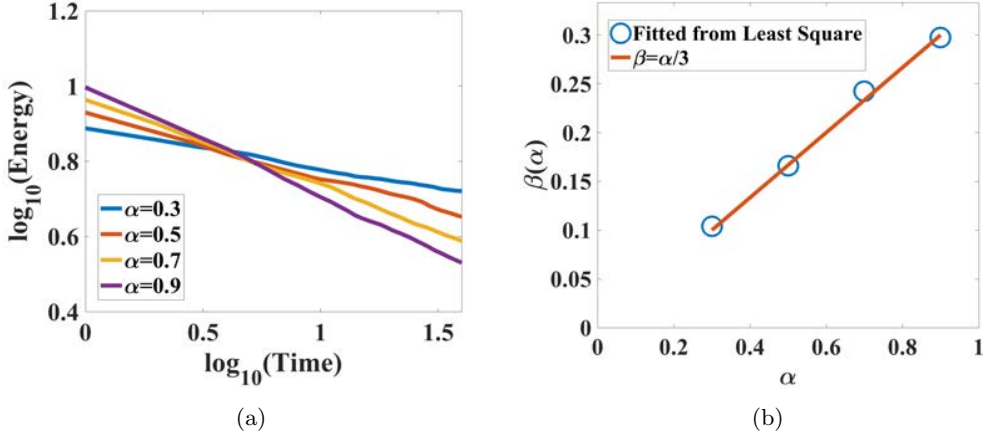

(a)                  (b)

Figure 4.5: Energy scaling laws for the time-fractional molecular beam epitaxy equation with slope selection. (A) The log-log plot of the effective energy and time with different fractional-order $\alpha$. (b) The least-square fitted energy decay power $\beta$ for different fractional order $\alpha$.

Overall, we lay out a foundation for time fractional phase field models by discovering a universal power laws scaling property. *The scaling of effective free energy/roughness in the time-fractional phase field models during coarsening follow a similar power law as the integer phase field models, where the power is linearly proportional to the fractional order $\alpha$.* This universal scaling law is verified numerically against several phase field models, including the Cahn-Hilliard (CH) equations with different variable mobilities and molecular beam epitaxy (MBE) models. This new finding sheds light on potential applications of time fractional phase field models in studying coarsening dynamics and crystal growths. A rigorous asymptotic analysis to this interesting correlation will be pursed in our future work.

# Acknowledgment

Lizhen Chen would like to acknowledge the support from National Science Foundation of China through Grant 11671166 and U1530401 , Postdoctoral Science Foundation of China through Grant 2015M580038. Jia Zhao is partially supported by a seed grant (Research Catalyst Grant) from Office of Research and Graduate Studies at Utah State University. Hong Wang is partially supported by the OSD/ARO MURI Grant W911NF-15-1-0562 and by the National Science Foundation under Grant DMS-1620194.

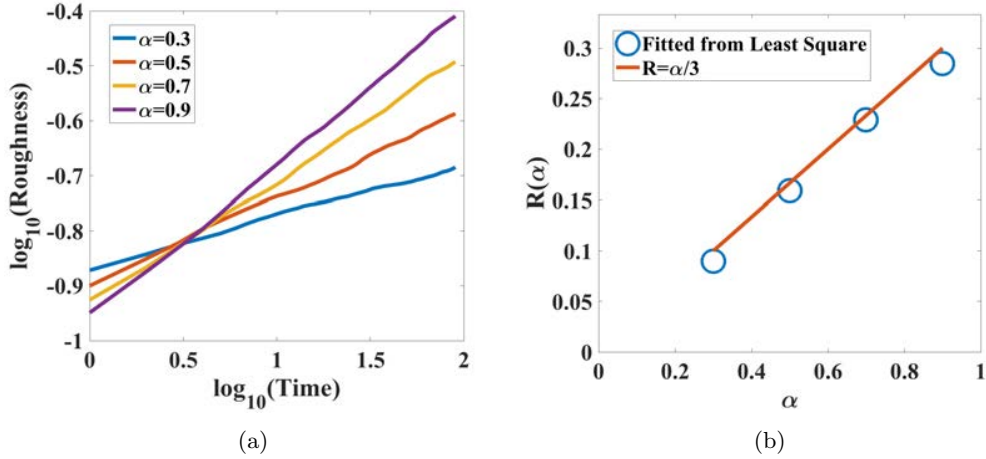

Figure 4.6: Roughness scaling laws for the time-fractional molecular beam epitaxy equation with slope selection. (A) The log-log plot of the roughness and time with different fractional-order $\alpha$. (b) The least-square fitted energy decay power $\beta$ for different fractional order $\alpha$.

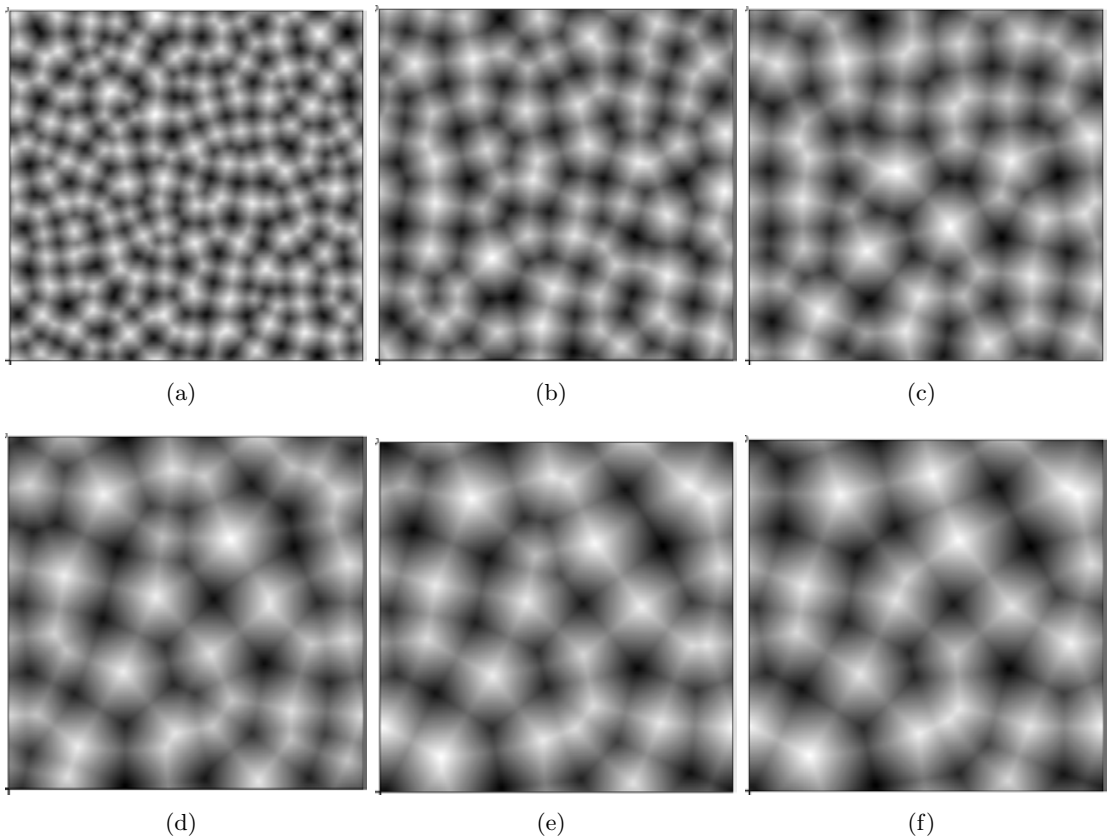

Figure 4.7: Time snapshots of the height evolution ($\phi$) of FMBE model with slope selection with $\alpha = 0.9$. The snapshots of $\phi$ at time $t = 5, 25, 50, 100, 150, 200$ are presented.

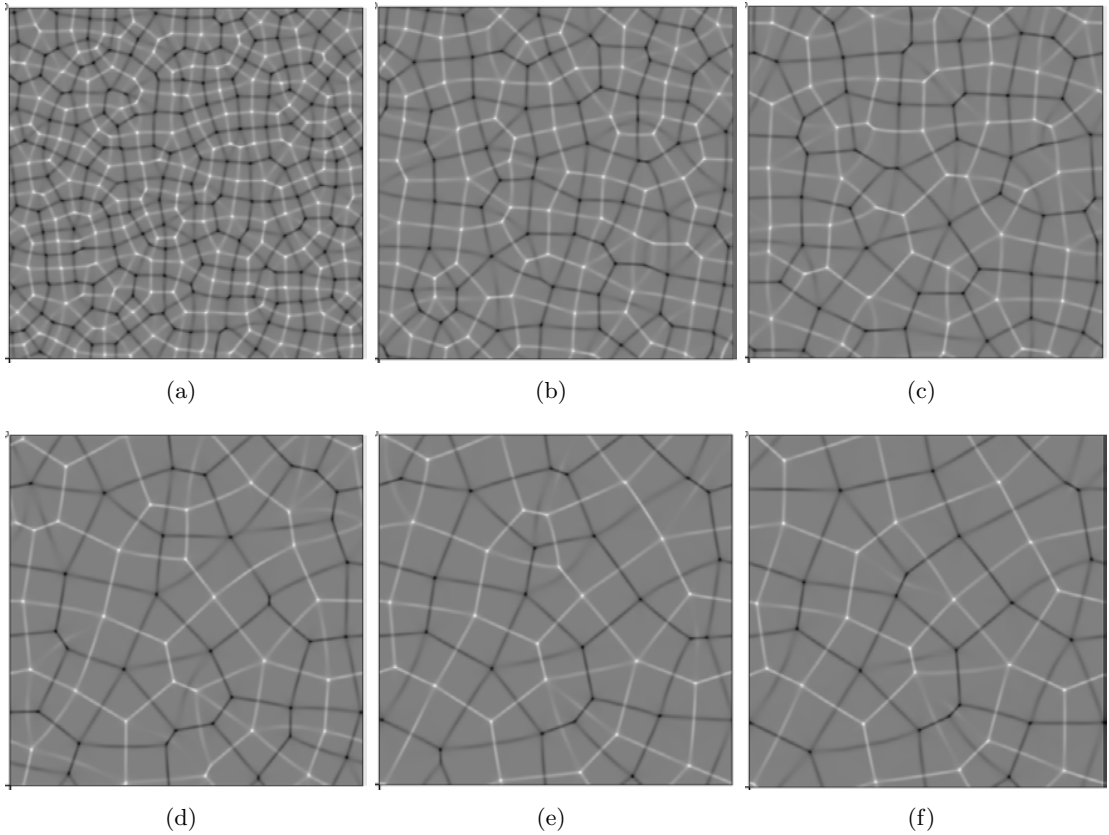

Figure 4.8: Time snapshots of the roughness evolution ($\Delta\phi$) of the FMBE model with slope selection with $\alpha = 0.9$. The snapshots of $\Delta\phi$ at time $t = 5, 25, 50, 100, 150, 200$ are presented.

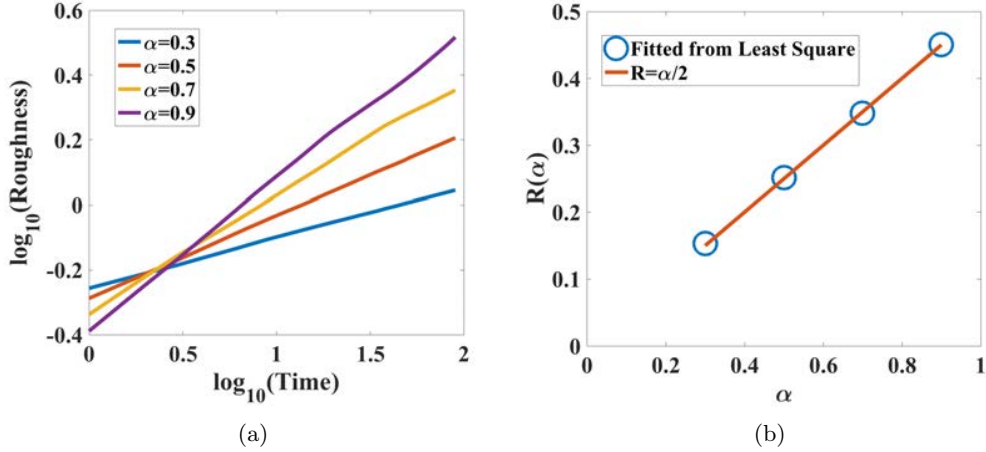

Figure 4.9: Roughness scaling laws for the FMBE model without slope selection. (A) The log-log plot of the roughness and time with different fractional-order $\alpha$. (b) The least-square fitted energy decay power $R(\alpha)$ for different fractional order $\alpha$.

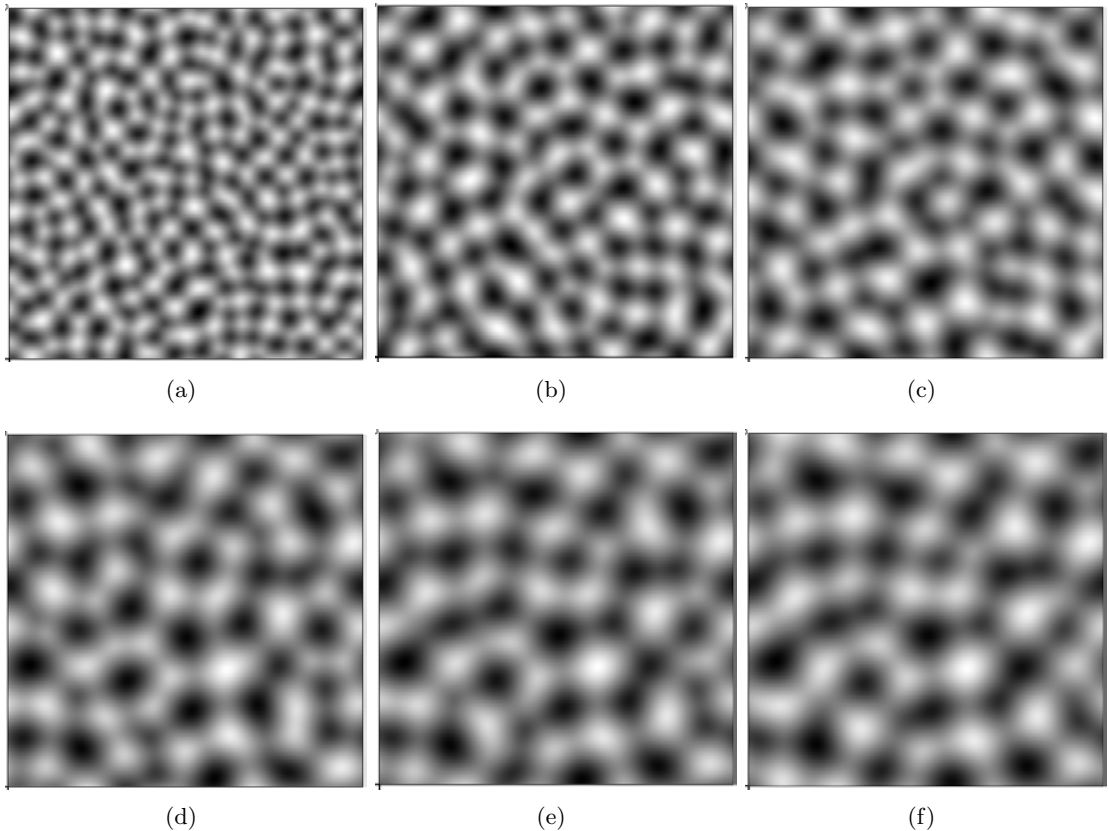

Figure 4.10: Time snapshots of the height evolution ($\phi$) of FMBE model without slope selection with $\alpha = 0.9$. The snapshots of $\phi$ at time $t = 5, 25, 50, 100, 150, 200$ are presented.

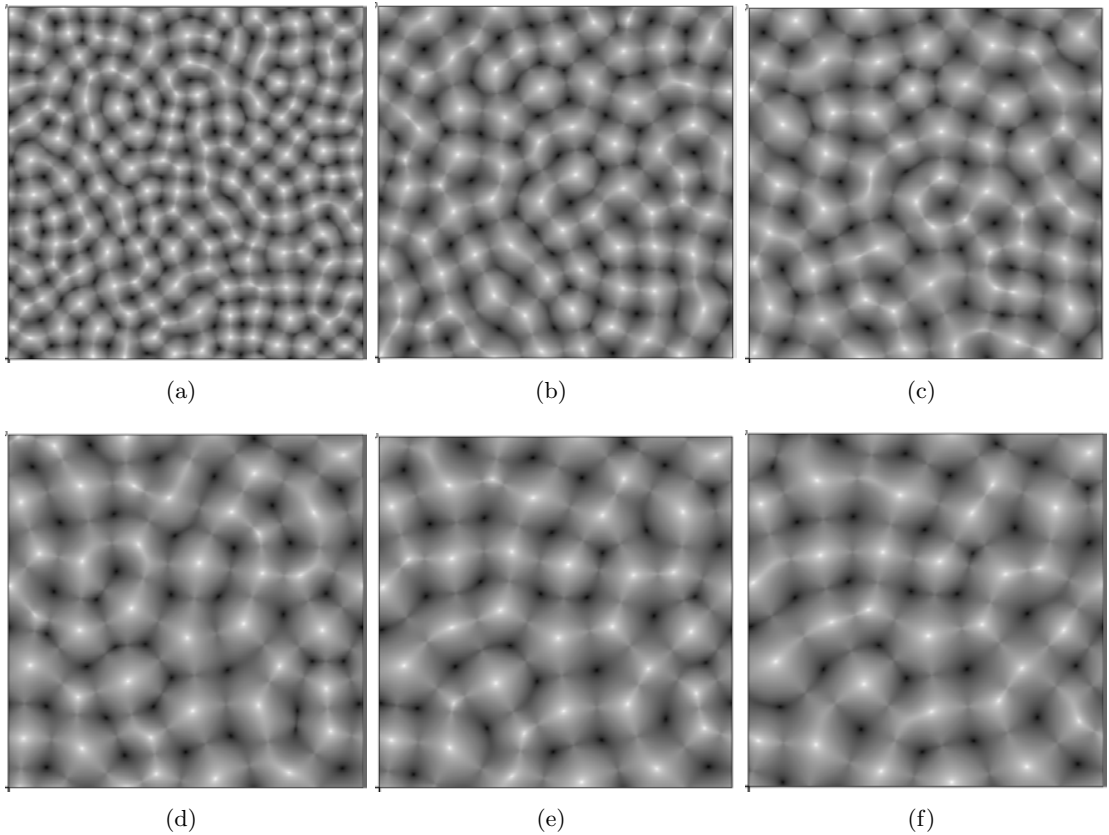

Figure 4.11: Time snapshots of the roughness evolution ($\Delta\phi$) of FMBE model without slope selection with $\alpha = 0.9$. The snapshots of $\Delta\phi$ at time $t = 5, 25, 50, 100, 150, 200$ are presented.